\documentclass[11pt]{amsart}
\usepackage{amsfonts}

\usepackage[dvips]{graphicx,psfrag}
\usepackage{multicol}

\newtheorem{theorem}{Theorem}[section]
\theoremstyle{plain}

\newtheorem{definition}{Definition}

\newtheorem{proposition}{Proposition}[section]
\newtheorem{remark}{Remark}[section]

\numberwithin{equation}{section}

\newcommand{\fimprova}{
\hspace*{\fill} \rule{0.15cm}{0.3cm} }

\begin{document}

\title[A Free Boundary Isoperimetric Problem in $\mathbb H^3$]{A Free Boundary Isoperimetric Problem in the Hyperbolic Space between parallel horospheres}

\author{Rosa M.B. Chaves}\address{Instituto de Matem\'atica e Estat\'\i stica \hfill\break\indent Universidade de S\~ao Paulo \hfill\break\indent Rua do Mat\~ao, 1010, S\~ao Paulo - SP, Brazil \hfill\break\indent CEP 05508-090} \email{rosab@ime.usp.br}
\author{Renato H.L. Pedrosa}\address{Instituto de Matem\'atica, Estat\'\i stica  e Computa\c c\~ao Cient\'ifica \hfill\break\indent Universidade Estadual de Campinas \hfill\break\indent Rua S\'ergio Buarque de Holanda, 651, Campinas - SP, Brazil \hfill\break\indent CEP: 13083-859} \email{pedrosa@ime.unicamp.br}
\author{M\'arcio F. da Silva}\address{Centro de Matem\'atica, Computa\c c\~ao e Cogni\c c\~ao \hfill\break\indent Universidade Federal do ABC \hfill\break\indent Rua Catequese, 242, Santo Andr\'e - SP, Brazil \hfill\break\indent CEP 09090-400} \email{marcio.silva@ufabc.edu.br}


\begin{abstract}
In this work  we investigate the following isoperimetric problem:
to find the regions of prescribed volume with minimal boundary
area between two parallel horospheres in hyperbolic 3-space (the
area of the part of the boundary contained in the horospheres is
not included). We reduce the problem to the study of rotationally
invariant regions and obtain  the  possible isoperimetric
solutions by studying  the behaviour of the profile curves of the
rotational surfaces with constant mean curvature in the hyperbolic
space.

\end{abstract}

\maketitle

\section{Introduction}

Geometric isoperimetric problems, (upper) estimates for the volume
of regions of a given fixed boundary volume, or the dual problems,
play an important role in Analysis and Geometry. There are both
isoperimetric inequalities , common in Analysis, and actual
classification of optimal geometric objects, like the round ball
in Euclidean Geometry. We will be interested in the study of a
relative free-boundary isoperimetric problem in Hyperbolic 3-space
between two parallel horospheres. A survey of recent results about
the geometric isoperimetric problems is \cite{RR}.

For a Riemannian manifold $M^{n}$, the classical isoperimetric
problem assumes the following formulation: to classify, up to
congruency by the isometry group of M, the (compact) regions
$\Omega \subseteq M$ enclosing a fixed volume that have minimal
boundary volume. The relevant concepts of volume involved are
those of Geometric Measure Theory: regions and their boundaries
are $n$-rectifiable (resp. $(n-1)$-rectifiable) subsets of $M$
(cf. \cite{Mo2}).

 If $M$ has boundary, the part of $\partial \Omega$
included in the interior of $M$ will be called the \emph{free
boundary} of $\Omega$, the other part will be called the
\emph{fixed boundary}. One may specify how the fixed boundary of
$\Omega$ is included in the computation of the boundary volume
functional. In this paper, we will assume that the volume of the
fixed boundary of $\Omega$ is not considered in the boundary
volume functional. We will see in Section 3 that this implies that
the angle of contact between the interior boundary of $\Omega$ and
$\partial M$ is $\pi/2$ (when such contact occurs). Such relative
problems are related to the geometry of stable drops in
capillarity problems (the angle of contact depends, as mentioned,
on how one considers the volume of the fixed boundary in the
computation of the boundary volume functional). For a discussion
about that, we refer to \cite{Finn}.

The motivation for our work is the well-known result of
Athanassenas \cite{At} and Vogel \cite{V} which implies that
between two parallel planes in the Euclidean space
$\mathbb{R}^{3}$, a (stable) soap-bubble touching both walls
perpendicularly is a straight cylinder perpendicular to the planes
(enclosing a tube), and may only exist down to a certain minimal
enclosing volume depending on the distance between the planes.
Below that value, only half-spheres touching one of the planes or
whole spheres not touching either plane occur, the cylinders
becoming unstable. A newer proof can be found in \cite{PR}, where
the authors study the analogous problem in higher-dimensional
Euclidean spaces.

In this paper we study the analogous relative isoperimetric
problem between two parallel horospheres in hyperbolic space
$\mathbb{H}^{3}(-1)$. We will use the upper halfspace model
$\mathbb{R}^{3}_{+}$. The parallel horospheres are then
represented by the horizontal Euclidean 2-planes of
$\mathbb{R}^{3}_{+}$.  We will present in this paper a detailed
classification of the possible isoperimetric solutions.

The existence of isoperimetric regions in the manifold with
boundary $(B,g)$, the slab composed by the two horospheres and the
region between them, may be obtained by adapting a result due to
Morgan \cite{Mo2} ($B/G$ is a compact space, where $G$ is the
subgroup of the isometry group of $\mathbb{H}^{3}(-1)$ leaving $B$
invariant, so that Morgan's result applies). Regarding the
regularity of the free boundary, well-known results about the
lower codimension bounds of the singular subset imply that it must
be regular, in fact analytic.

In Section 2 we give some basic definitions in the model
$\mathbb{R}^{3}_{+}$ like geodesics, totally geodesic surfaces,
umbilical surfaces and  rotational surfaces. We also give a more
precise formulation for the isoperimetric problem considered using
area and volume functionals.

In Section 3 we discuss briefly the rotationally invariance of
isoperimetric regions and some of their basic geometric
properties, since their (free) boundaries must have constant mean
curvature and, when touching the bounding horospheres, the contact
angle must be $\pi/2$. We will also discuss in more detail the
existence of isoperimetric regions and the regularity of the free
boundary part.

In Section 4 we investigate the tangency of profile curves for the
rotational surfaces with constant mean curvature to determine the
possible isoperimetric regions between the two parallel
horospheres in $\mathbb{R}^{3}_{+}$. We define what we mean by
catenoid, equidistant and onduloid type surfaces and prove the
following result.

Let $c_{1},  c_{2}$ be positive real constants, $c_{1} < c_{2}$,
and $\mathcal{F}_{c_{1}, c_{2}} = \{(x,y,z) \in
\mathbb{R}^{3}_{+}: c_{1} \leq z \leq c_{2} \}$. Let $V>0$ and
$\mathcal{C}_{c_{1},c_{2},V}$ be the set of $ \Omega \subset
\mathcal{F}_{c_{1}, c_{2}}$ with volume $|\Omega|=V$ and boundary
volume (area) $A(\Omega \bigcap \stackrel{o}{\mathcal{F}_{c_{1},
c_{2}}}) < \infty $, where we suppose $\Omega$ to be connected,
compact, 3-rectifiable in $\mathcal{F}_{c_{1}, c_{2}}$, having as
boundary (between the horospheres) an embedded, orientable,
2-rectifiable surface.

\begin{theorem} Let $A_{c_{1},c_{2},V} = \inf \{A(\Omega \bigcap
\stackrel{o}{\mathcal{F}_{c_{1},c_{2}}}): \Omega \in
\mathcal{C}_{c_{1},c_{2},V} \}$. Then
\begin{enumerate}
    \item there exists $\Omega \in \mathcal{C}_{c_{1},c_{2},V}$
    such that $A(\Omega \bigcap \stackrel{o}{\mathcal{F}_{c_{1},c_{2}}}) =
    A_{c_{1},c_{2},V}$. As already mentioned, the free boundaries are actually analytic
    surfaces;
    \item if $\Omega$ has minimal boundary volume, between the horospheres, the free boundary of $\Omega$, is either
    \begin{enumerate}
            \item of catenoid type or umbilical with $H=1$, or
        \item of equidistant type or umbilical with $0 < H <1$, or
        \item of onduloid type or umbilical with $H > 1$.
        \end{enumerate}
\end{enumerate}
\end{theorem}

The details of the description above are included in Section 4. We
observe here that this result shows how the situation in
hyperbolic geometry is different from the one in the Euclidean
3-space, where we also have rotationally invariant surfaces of
catenoid and of onduloid type, but they cannot appear as
boundaries of optimizing tubes (even though, in higher dimensions,
hypersurfaces generated by onduloids in Euclidean spaces are known
to occur as boundaries of optimal tubes connecting two parallel
hyperplanes (cf. \cite{PR})).

\noindent \textbf{Acknowledgements}. The author M\'arcio Silva
would like to thank CAPES and CNPq for support during his PhD
studies, which included part of the work presented in this paper.

\section{Preliminaries}
In this section we will introduce some basic facts and notations
that will appear along the paper.

Let $\mathcal{L}^4=(\mathbb{R}^4,g)$ the $4$-dimensional Lorentz
space endowed with the metric $g(x,y)=\ x_1 y_1 + x_2 y_2 + x_3
y_3 - x_4 y_4$ and the $3$-dimensional hyperbolic space
\begin{displaymath}
\mathcal{H}^{3}(-1):=\{p=(x_1,x_2,x_3,x_4)\in \mathcal{L}^4 : \
g(p,p)=-1,\ x_4>0\}.
\end{displaymath}
We use the upper halfspace model $\mathbb{R}^{3}_{+}:= \{(x,y,z)
\in \mathbb{R}^{3}; z>0 \}$ for $\mathcal{H}^{3}(-1)$, endowed
with the metric
$<,>=ds^{2}=\displaystyle\frac{dx^{2}+dy^{2}+dz^{2}}{z^{2}}.$

Let $\phi:\Sigma \rightarrow \mathbb{R}^{3}_{+}$ be an isometric
immersion of a compact manifold $\Sigma$ with boundary
$\partial\Sigma \neq \emptyset$ and  $\Gamma$ be  a curve in
$\mathbb{R}^{3}_{+}$. If $\phi$ is a diffeomorphism of
$\partial\Sigma$ onto $\Gamma$, we say that  $\Gamma$ is the
boundary of $\phi$ and if $\phi$ has constant mean curvature $H$,
we say that $\Sigma$ is an $H$-surface with boundary $\Gamma$. We
identify $\Sigma$ with its image by $\phi$ and $\partial\Sigma$
with the curve $\Gamma$.

The plane $z=0$ is called the infinity boundary of
$\mathbb{R}^{3}_{+}$ and we denote it  by
$\partial_{\infty}\mathbb{R}^{3}_{+}$. The geodesics of
$\mathbb{R}^{3}_{+}$ are represented by vertical Euclidean lines
and half-circles orthogonal to
$\partial_{\infty}\mathbb{R}^{3}_{+}$, contained in
$\mathbb{R}^{3}_{+}$. The totally geodesic surfaces have constant
mean curvature $H=0$ and are represented by vertical Euclidean
planes and hemispheres orthogonal to
$\partial_{\infty}\mathbb{R}^{3}_{+}$, contained in
$\mathbb{R}^{3}_{+}$.

The horizontal Euclidean translations and the rotations around a
vertical geodesic are isometries of $\mathbb{R}^{3}_{+}$. We have
two family of isometries associated to one point $p_{0}\in
\partial_{\infty}\mathbb{R}^{3}_{+}$: the Euclidean homotheties centered in $p_{0}$ with factor $k > 0$, called
hyperbolic translations through a geodesic $\alpha$ perpendicular
to $\partial_{\infty}\mathbb{R}^{3}_{+}$ in $p_{0}$, and the
hyperbolic reflections with respect to a totally geodesic surface
$P$.

When $P$ is a hemisphere orthogonal to
$\partial_{\infty}\mathbb{R}^{3}_{+}$ centered in $p_{0}$ and
radius $r > 0$, the hyperbolic reflections are Euclidean
inversions centered in $p_{0}$ that fix $P$ and, when $P$ is a
vertical Euclidean plane, they are Euclidean reflections with
respect to $P$.

The umbilical surfaces of $\mathbb{R}^{3}_{+}$ are described as
follows (see for example \cite{L}):
\begin{enumerate}
\item {\it Totally geodesics:} These surfaces were already
described ($H=0$).

\item {\it Geodesic spheres:} represented by Euclidean spheres
entirely contained in $\mathbb{R}^{3}_{+}$, they have $H>1$ (mean
curvature vector points to the interior). If $\rho$ is the
hyperbolic radius of a geodesic sphere then  $H=\coth \rho$.

\item {\it Horospheres:} represented by horizontal Euclidean
planes of $\mathbb{R}^{3}_{+}$ and Euclidean spheres of
$\mathbb{R}^{3}_{+}$ which are tangent to
$\partial_{\infty}\mathbb{R}^{3}_{+}$, they have $H=1$ and the
mean curvature vector points upward in the case of horizontal
planes and to the interior in the  case of spheres.

\item {\it Equidistant surfaces:} represented by the intersection
of $\mathbb{R}^{3}_{+}$ with the planes of $\mathbb{R}^{3}$ that
are neither parallel nor perpendicular to the plane ${z=0}$ and by
(pieces of) Euclidean spheres that are not entirely contained in
$\mathbb{R}^{3}_{+}$ and are neither parallel tangent nor
perpendicular to the plane ${z=0}$. They have $0 < H < 1$ and the
mean curvature vector points to the totally geodesic surface they
are equidistant to.

\end{enumerate}

In our study, the (spherical) rotational surfaces of
$\mathbb{R}^{3}_{+}$ play an important role since the solutions of
the isoperimetric problem must be rotationally invariant regions.
 They are defined as surfaces invariant by a subgroup of isometries whose
principal orbits are (Euclidean) circles.

Let $\Pi_{1}$ and  $\Pi_{2}$ be horospheres represented by
distinct parallel horizontal Euclidean planes,
  $\Pi = \Pi_{1} \cup \Pi_{2}$, $\mathcal{F}=\mathcal{F}(\Pi_1,\Pi_2)$ the closed slab between them and $\phi: \Sigma \rightarrow \mathcal{F}$
an immersion of a compact, connected, embedded and  orientable
$C^{2}$- surface with boundary $\Gamma =
\partial \Sigma$ and such that $\phi(\Gamma)\subset \Pi$ (we will see later that the image by
$\phi$ of the interior of the surface $\Sigma$ will not be allowed
to touch $\Pi$ if $\Sigma$ is the boundary of an optimal domain in
our variational problem, but this is not part of the general
situation at this point).

We fix, now, the notation for some well-known geometric invariants
related to isometric immersions. We identify (locally) $\Sigma$
with $\phi( \Sigma)$ and  $X(p) \in T_{p} \Sigma$ with $d
\phi_{p}(X(p)) \subset \mathbb{R}^{3}_{+}$. We have, as usual, the
decomposition $T_{p}(\mathbb{R}^{3}_{+}) = T_{p}(\Sigma) \oplus
N_{p}(\Sigma)$ into the tangent and normal spaces to $\Sigma$ in
$p$, respectively. Choose an orientation for $\Sigma$ and let $N$
be the (positive) unitary normal field along the immersion $\phi$.
If $X(p) \in T_{p}(\mathbb{R}^{3}_{+})$, we may write $X(p) =
X(p)^{T} + X(p)^{N}=X(p)^{T} + \alpha N(p),$ where $\alpha \in
\mathbb{R}$.

Let $<\!\cdot,\cdot\!>$ be the metric  induced on $\Sigma$ by the
immersion $\phi$, making it isometric, $\overline{\nabla}$ be the
Riemannian connection of the ambient space $\mathbb{R}^{3}_{+}$
and $\nabla$ the induced Riemannian connection of  $\Sigma$. Let
$X, Y \in \mathfrak{X}(\Sigma)$ be $C^\infty$-vector fields, then
$\nabla_{X} Y = (\overline{\nabla}_{X} Y)^{T}$ and $\mathcal{B}(X,
Y) = (\overline{\nabla}_{X} Y)^{N}$ are, as usual, the
\emph{induced connection} on $\Sigma$ and the \emph{second
fundamental form} of the immersion given by $\mathcal{B}$. We also
have the \emph{Weingarten operator}
$A_{N}(Y)=-(\overline{\nabla}_{Y} N)^{T}$ so that
$<\!A_{N}(X),Y\!>=<\!\mathcal{B}(X,Y),N\!>.$ Finally, the mean
curvature of the immersion $\phi$ is $H  =1/2 ~
\textrm{trace}(A_{N})$.

\begin{definition}
A variation of $\phi$ is a smooth map $F: (-\epsilon, \epsilon)
\times \Sigma \rightarrow \mathbb{R}^{3}_{+}$, such that for all
$t \in (-\epsilon, \epsilon)$, $\phi_{t}: \Sigma \rightarrow
\mathbb{R}^{3}_{+}$,  defined by
$\phi_{t}(p) = F(t,p)$,\\
(a)~ is an immersion;\\
(b)~ $\phi_{0} = \phi$.
\end{definition}

For $p \in \Sigma$, $X(p) = \partial\phi_{t}(p)/\partial
{t}|_{t=0}$ is the \emph{variation vector field} of $F$ and the
\emph{normal variation function} of $F$ is given by $f(p) =
<\!X(p),N(p)\!>$. We say that the variation $F$ is {\it normal} if
$X$ is normal to $\phi$ at each point and that $F$ has compact
support if
 $X$ has compact support. For a variation with compact support and small values of $t$ we have that  $\phi_{t}$ is an immersion of $\Sigma$ in
$\mathbb{R}^{3}_{+}$. In this case we define  the {\it (variation)
area  function} $A: (-\epsilon, \epsilon) \rightarrow \mathbb{R}$
by
\begin{displaymath}
A(t) = \displaystyle\int_{\Sigma} dA_{t} =
\displaystyle\int_{\Sigma}\sqrt{\det[
(d\phi_{t})^{\ast}(d\phi_{t})]}~ dA,
\end{displaymath}
where $dA$ is the element of area of $\Sigma$. $A(t)$ is  the area
of $\Sigma$ with the metric induced by $\phi_{t}$. We also define
the {\it (variation) volume function} $V: (-\epsilon, \epsilon)
\rightarrow \mathbb{R}$ by
\begin{displaymath}
V(t) = -\displaystyle\int_{[0,t] \times \Sigma} F^{\ast}
d(\mathbb{R}^{3}_{+}),
\end{displaymath}
where  $d(\mathbb{R}^{3}_{+})$ is the element of  volume of
$\mathbb{R}^{3}_{+}$ and $F^{\ast} d(\mathbb{R}^{3}_{+})$ is the
pull-back of $d(\mathbb{R}^{3}_{+})$ by $F$. $V(t)$ is not,
actually, the volume of some region with $\phi_t(\Sigma)$ as
boundary, but of a ``tubular neighborhood'' along $\phi(\Sigma)$
between $\phi(\Sigma)$ and $\phi_t(\Sigma)$. The sign is related
to the net change with respect to the normal field defining the
orientation (for example, contracting a sphere in $\mathbb{R}^3$,
which means moving it in the direction of the mean curvature
vector, gives a negative sign for $V(t)$, as expected).

\begin{definition} \label{def:variation}
Let $F: (-\epsilon, \epsilon) \times \Sigma \rightarrow
\mathbb{R}^{3}_{+}$ a variation of $\phi$.\\
$(i)$ $F$ preserves volume if $V(t) = V(0) (=0), \forall t \in
(-\epsilon, \epsilon)$;\\
$(ii)$ $F$ is admissible  if $F(\partial \Sigma) \subset \Pi,
\forall t \in (-\epsilon, \epsilon)$.
\end{definition}

\begin{definition} \label{def:stationary}
We say that the immersion $\phi$ is stationary when $A'(0)=0$, for
all admissible variations that preserve volume.
\end{definition}

\begin{remark} If $\Omega$ is a (compact) regular region in the slab $\mathcal{F}$ between the horospheres $\Pi$, taking $\Sigma$ in
Definition \ref{def:variation} as the (embedded regular)  free
boundary of $\Omega$, we may extend the above variational approach
to produce a variation $\Omega(t)$ of $\Omega$ by embedded domains
(for small $t$), such that the condition $V(t)=0$ in Definition
\ref{def:variation} is equivalent to holding the measure (volume)
$|\Omega(t)|$ of $\Omega(t)$ constant (the same as that of
$\Omega=\Omega(0)$) along the variation. This justifies saying
that the variation ``preserves volume'' in the above definition.
\end{remark}

We end this section by stating again our problem.  Let $\Pi_1$ and $\Pi_2$ be two parallel horospheres in $\mathbb{R}^{3}$ and $\mathcal{F}=\mathcal{F}(\Pi_2,\Pi_2)$ the (closed) slab between them. \\

\noindent \textbf{Isoperimetric problem for
$\mathcal{F}(\Pi_1,\Pi_2)$}: \emph{fix a volume value and study
the domains $\Omega\subset \mathcal{F}$ with the prescribed volume
which have minimal free boundary area.}

\begin{definition}\label{def:isopdomain}
A (compact) minimizing region $\Omega$ for this problem will be
called an \textbf{isoperimetric domain} or \textbf{region} in
$\mathcal{F}$.
\end{definition}

In more detail, one is looking for \emph{the classification and
geometric description of isoperimetric regions (as a function of
the volume value), in as much detail as possible, aiming at the
determination of
the isoperimetric profile (minimal free boundary area as a function of prescribed volume) for $\mathcal{F}$.}\\



\section{First results about the isoperimetric solutions}

Our main goal in this section is to characterize the stationary
immersions according to Definition \ref{def:stationary}. The
following formulae for the first variations of the area and volume
functions are well known. For an immersed surface with boundary,
the \emph{exterior conormal} is the vector field along the
boundary given as follows: in the tangent plane of  $\Sigma$ in $p
\in \Gamma$, take the perpendicular vector to the tangent vector
to $\Gamma$ in $p$  along the boundary $\Gamma$.

\begin{proposition}
\label{prop:firstvariations} Let $F$ be a variation of $\phi$ with
variational field $X$ and compact support in $\Sigma$. Then
\begin{enumerate}
\item $A'(0)= -2 \displaystyle\int_{\Sigma}H  f ~ dA +
\displaystyle\int_{\Gamma} <X , \nu > d\Gamma$, where $\nu$ is the
unitary exterior conormal, $dA$ is the element of area of $\Sigma$
and $d\Gamma$ is the element of length of $\Gamma$ induced by
$\phi$.
\item $ V'(0)= -  \displaystyle\int_{\Sigma} f~ dA,$ where  $f(p)
=~ <\!X(p),N(p)\!>$, as before.
\end{enumerate}
\end{proposition}
\textbf{Proof}:~ Although the formula of the variation of the area
functional is well known (see \cite {BCE}), we show here a
different proof to get it. From the definition of $A(t)$ we obtain
\begin{displaymath}  \begin{array}{rl}
A'(t) = &{\displaystyle \int_{\Sigma}} \left[ \frac{1}{2
\sqrt{\det [(d\phi_{t})^{\ast} d\phi_{t}]}} \Big( \det
[(d\phi_{t})^{\ast} d\phi_{t}] \Big). \right. \\
& \left. \textrm{trace} \Big( [(d\phi_{t})^{\ast} d\phi_{t}]^{-1}
\circ \frac{d}{dt} ((d\phi_{t})^{\ast} d\phi_{t}) \Big)\right]~
dA.
\end{array}
\end{displaymath}

As $\phi_{0}$ is the inclusion of $\Sigma$ in
$\mathbb{R}^{3}_{+}$, $d\phi_{0}$ is the inclusion of the
respective tangent spaces and $d\phi_{0}^{\ast}$ is the orthogonal
projection on $T\Sigma$.

Evaluating $A'(t)$ for $t=0$ we get
\begin{displaymath}
A'(0) = \displaystyle\int_{\Sigma} \displaystyle\frac{1}{2}
\textrm{trace} \Big( \displaystyle\frac{d}{dt}\Big|_{t=0}
[(d\phi_{t})^{\ast} d\phi_{t}] \Big) dA.
\end{displaymath}

Using the Symmetry Lemma for the connection $\nabla^{\phi}$ along
the immersion we have
\begin{displaymath}
\displaystyle\frac{d}{dt}\Big|_{t=0} (d\phi_{t}) = \nabla^{\phi}
\displaystyle\frac{\partial \phi_{t}}{\partial{t}} \Big|_{t=0} =
\nabla^{\phi} X.
\end{displaymath}

Then
\begin{displaymath}
A'(0) = \displaystyle\int_{\Sigma} \displaystyle\frac{1}{2}
\textrm{trace} \Big( (\nabla^{\phi} X)^{\ast} \Big|_{T\Sigma} +
\textrm{proj}_{T\Sigma} \nabla^{\phi} X \Big)
dA\\
=\displaystyle\int_{\Sigma} \textrm{trace} \Big(
\textrm{proj}_{T\Sigma} \nabla^{\phi} X \Big) dA,
\end{displaymath}
where $\textrm{proj}_{T\Sigma}$ denotes the projection on
$T\Sigma$.

Decomposing the variational field as $X = X^{T} + X^{N}$, we have
that the projections of the tangent and normal components of
$\nabla^{\phi}( X )$ on $T\Sigma$ are
\begin{displaymath}
\textrm{proj}_{T\Sigma} \nabla^{\phi}( X^{T}) = \nabla ( X^{T})
\end{displaymath}
\begin{displaymath}
\textrm{proj}_{T\Sigma}\nabla^{\phi}( X^{N}) = - A_{X^{N}}.
\end{displaymath}

So
\begin{displaymath}
A'(0) = \displaystyle\int_{\Sigma} \Big( div X^{T} - 2 <\!X^{N} ,
H~ N\!> \Big) dA.
\end{displaymath}

Using the Stokes theorem we get
\begin{displaymath}
A'(0)= \displaystyle\int_{\Gamma} <\!X^{T} , \nu \!
> d\Gamma -2 \displaystyle\int_{\Sigma}<\!X^{N} , H~ N\!> dA =
\end{displaymath}
\begin{displaymath}
=~ -2 \displaystyle\int_{\Sigma} H f~ dA +
\displaystyle\int_{\Gamma} <\!X , \nu \!> d\Gamma.
\end{displaymath}

The first variation of volume given in (2) is the standard one and
will be omitted (cf.\cite{BCE}).

\fimprova

 From the  next result we  deduce that the boundary of the
isoperimetric region we are studying must be an $H$-surface that
makes a contact angle $\pi/2$ with the horospheres $\Pi_1$ and
$\Pi_2$.

\begin{theorem}
An immersion $\phi: \Sigma \rightarrow \mathbb{R}^{3}_{+}$ is
stationary if and only if $\phi$ has constant mean curvature and
meets the horospheres $\Pi =\Pi_1 \cup \Pi_2$ that contains its
boundary $\Gamma =
\partial \Sigma$, perpendicularly along the boundary (if the intersection is non-empty).
\end{theorem}
\textbf{Proof}: Adapting the proof of Proposition 2.7 of \cite{BC}
we show that if $\phi$ has constant mean curvature and meets the
horospheres $\Pi =\Pi_1 \cup \Pi_2$ that contains its boundary
$\Gamma = \partial \Sigma$, perpendicularly along the boundary
then  $\phi$ is stationary; and that if $\phi$ is stationary then
$\phi$ has constant mean curvature. To show that if $\phi$ is
stationary then $\phi$ meets $\Pi$ perpendicularly along its
boundary $\Gamma =
\partial \Sigma$, we take an admissible variation $\Phi$ that
preserves volume with variational field
 $X$ and $p_{0} \in \partial \Sigma$. Suppose, by contradiction,
$<\! X(p_{0}),\nu(p_{0})\!> \neq 0$. By continuity there is a
neighbourhood
 $U =  W_{1} \bigcap \partial \Sigma $ of $p_{0}$ such that
$<\!X(p),\nu(p)\!> ~ > 0$, $\forall p \in U$, where $W_{1}$ is a
neighbourhood of  $p_{0}$ in $\Sigma$. We take $q \in ~
\stackrel{o}{\Sigma} \backslash W_{1}$ , $W_{2}$ a neighbourhood
of $q$ disjoint to $W_{1}$ and a partition of unity on $W_{1}
\bigcup W_{2}$. There exists a differentiable function
$\xi_{1}:W_{1} \rightarrow \mathbb{R}$ such that $\xi_{1}(W_{1})
\subset [0,1]$ with support $\textrm{supp} ~ \xi_{1} \subset
W_{1}$. We may consider also a differentiable map
$\xi_{2}:W_{2}\rightarrow \mathbb{R}$  such that $\xi_{2}(W_{2})
\subset
 [0,1],~
 \textrm{supp} ~ \xi_{2} \subset W_{2}$ and
 \begin{displaymath}\label{perpendbordo}
    \displaystyle\int_{W_{1}}\xi_{1}f ~ dW_{1}+ \displaystyle\int_{W_{2}}\xi_{2}f ~
    dW_{2}=0.
\end{displaymath}
We consider the variation $\Phi_{\xi}: (-\epsilon, \epsilon)
\times \Sigma \rightarrow \mathbb{R}$ with compact support on
$W_{1} \bigcup W_{2}$ such that
\begin{displaymath}
\Phi_{\xi}^{t}(p) = \Phi_{\xi}(t,p)=\left\{ \begin{array}{ll}
                         \Phi(\xi_{1} ~ t,p), ~ p \in W_{1},\\
                         \Phi(\xi_{2} ~ t,p), ~ p \in W_{2}.
                         \end{array}
                 \right.
\end{displaymath}

Notice that $\Phi_{\xi}$ is admissible because $\Phi$ is
admissible.

If $f_{\xi}(p)$ denotes the normal component of the variation
vector we have
\begin{displaymath}\label{perpmedianula}
    \displaystyle\int_{\Sigma} f_{\xi} ~ dA =
\displaystyle\int_{W_{1}} \xi_{1}(p) f(p) ~ dW_{1} +
\displaystyle\int_{W_{2}} \xi_{2}(p) f(p) ~ dW_{2} = 0,
\end{displaymath}
and $\Phi_{\xi}$ preserves volume.

For this variation we have
\begin{displaymath}\label{contraperp}
    0 = A'(0) = -2 H \displaystyle\int_{\Sigma} f_{\xi} ~ dA +
\displaystyle\int_{W_{1}} \xi_{1} ~ <\!X , \nu \!> d\Gamma =
\displaystyle\int_{W_{1}} \xi_{1} ~ <\!X , \nu \!> d\Gamma >0,
\end{displaymath}
which is a contradiction. Then $\forall p \in
\partial\Sigma$ it follows that $<\!X, \nu\!>(p) =0.$
\fimprova

 Next we show that the boundary of the isoperimetric solutions
are rotationally invariant regions.

We need some symmetrization principle for $H$-surfaces. Taking the
version of Alexandrov's Principle of Reflection for  the
hyperbolic space (for further references and details see \cite{A})
and using \cite{EB} as reference to take the suitable objects in
our case as the reflection planes, we get the next result (a
detailed proof may be found in \cite{L2}).

\begin{theorem}\label{teo:Alex}
Let $\Sigma$ be a compact connected orientable and  embedded
$H$-surface of class $C^{2}$, between two parallel horospheres
$\Pi_1 , \Pi_2$ in $\mathbb{R}^{3}_{+}$ and with boundary
$\partial \Sigma \subset \Pi_1 \bigcup \Pi_2$ (possibly empty).
Then $\Sigma$ is rotationally symmetric around an axis
perpendicular to $\Pi_1$ and $\Pi_2$.
\end{theorem}

We observe that the intersection of $\Sigma$ with a horosphere
$\mathcal{H}$ (represented by a horizontal Euclidean plane) is
just an Euclidean circle. In fact, if they were two concentric
circles and the isoperimetric region was delimited by these
circles we would
 get a totally geodesic symmetry plane $P$ that would not contain the axis of symmetry.

\section{Isoperimetric regions between horospheres in $\mathbb{R}^{3}_{+}$}

In this section we classify the rotational $H$-surfaces of
$\mathbb{R}^{3}_{+}$ between two parallel horospheres with
boundary contained in the horospheres and that intersects  the
horospheres perpendicularly. Therefore we get the possible
solutions for the isoperimetric problem in the hyperbolic space
since they must be regions delimited by those surfaces. We start
with some important results obtained by Barrientos \cite{B} in his
PhD thesis which will be useful for our task.

If $(\rho, \theta, z)$ are the cylindrical coordinates of a point
$p$ in $\mathbb{R}^{3}_{+}$ then the cartesian coordinates are
given by
\begin{equation}\label{cartescilin}
   (\overline{x},\overline{y},\overline{z})=e^{z}(\tanh \rho  \cos \theta, \tanh \rho  \sin \theta, \textrm{sech}
   \rho).
\end{equation}

For a spherical rotational surface $\Sigma$ of
$\mathbb{R}^{3}_{+}$ we can provide the so called  natural
parametrization, whose  metric is
\begin{displaymath}
    d\sigma^{2}=ds^{2}+U^{2}(s) ~ dt^{2},
\end{displaymath}
where $U=U(s)$ is a positive function. Supposing that its profile
curve is locally a  graphic $z= \lambda = \lambda (\rho)$ in the
plane $\theta =0$, the natural parameters are given by
\begin{displaymath}
    ds=\sqrt{1+\dot{\lambda}^{2}(\rho) ~ \cosh^{2}\rho} ~~
    d\rho ~~~ e ~~~ dt=d\varphi,
\end{displaymath}
and the following relations hold
\begin{equation}\label{coordnat}
    U^{2}(s)=\sinh^{2}\rho(s) ~~~ and ~~~
    \dot{\lambda}^{2}(s)=\displaystyle\frac{1+U^{2}(s)-\dot{U}^{2}(s)}{(1+U^{2}(s))^{2}}.
\end{equation}

Then  the natural parametrization  for a rotational surface in
cylindrical coordinates is
\begin{displaymath}
\left\{\begin{array}{l}
\sinh^{2}\rho(s) = U^{2}(s),\\
~\\
\lambda(s)=\displaystyle\int_{0}^{s}\frac{\sqrt{1+U^{2}(t)-\dot{U}^{2}(t)}}{1+U^{2}(t)} ~ dt,\\
~\\
\varphi(t)=t.
\end{array}
\right.
\end{displaymath}

In \cite{B}, the rotational were studied $H$-surfaces of
$\mathbb{R}^{3}_{+}$. By replacing $z(s)=U^2(s)$, he got the
differential equation for a rotational $H$-surface in
$\mathbb{R}^{3}_{+}$
\begin{displaymath}
\displaystyle\frac{\dot{z}^{2}}{4}=(1-H^{2})z^{2}+(1+2aH)z-a^{2},
\end{displaymath}
where $a$ is a constant of integration, and showed that the
behaviour of their profile curves is determined by the constant
$a$. Choosing  the orientation for the surfaces in order to $H
\geq 0$, there are three cases to study: $H=1$, $0\leq H <1$ and
$H>1$.

The natural parametrization for a rotational $H$-surface in
$\mathbb{R}^{3}_{+}$, ~$H=1$, generated by a curve
$c(s)=(\rho(s),\lambda(s))$ is given by
\begin{equation}\label{paramH1}
\left\{\begin{array}{lll}
\!\!\!\sinh^{2}\rho(s) = \displaystyle\frac{a^{2}+(1+2a)^{2}s^{2}}{1+2a},\\
~\\
\!\!\!\lambda(s)\!\!
=\!\!\displaystyle\int_{0}^{s}\!\frac{\sqrt{1+2a}(-a(1+a)+(1+2a)^{2}
t^{2}) \sqrt{a^{2}+(1+2a)^{2} t^{2}}}
{(-a(1+a)+(1+2a)^{2} t^{2})^{2}+(1+2a)^{4}t^{2}}~ dt,\\
~\\
\!\!\! \varphi(t) =t
\end{array}
\right.
\end{equation}

The natural parametrization for a rotational $H$-surface in
$\mathbb{R}^{3}_{+}$, ~ $0 \leq H <1$, generated by a curve
$c(s)=(\rho(s),\lambda(s))$ is given by
\begin{equation}\label{paramHrotaequid}
\left\{\begin{array}{lll}
\!\!\!\sinh^{2}\rho(s) = \displaystyle\frac{-A+B\cosh(2 \alpha s)}{2 \alpha^{2}},\\
~\\
\!\!\!
\lambda(s)\!=\!\!\!\displaystyle\int_{0}^{s}\!\frac{\sqrt{2}
\alpha (-2a\!\!+\!\!H(-1\!+\!B \cosh(2 \alpha t)))
\sqrt{\!-A\!\!+\!\!B \cosh(2 \alpha t)}}{(-2a\!\!+\!\!H (-1\!\!+\!\!B \cosh(2 \alpha t)))^{2}\!\!+\!\! \alpha^{2} B^{2} \sinh^{2}(2 \alpha t)} ~ \!dt,\\
~\\
\!\!\! \varphi(t)=t,
\end{array}
\right.
\end{equation}
where  $A=1+2aH, ~ B=\sqrt{1+4aH+4 a^{2}}$ and
$\alpha=\sqrt{1-H^{2}}.$

The natural parametrization for a rotational $H$-surface in
$\mathbb{R}^{3}_{+}$, ~ $H>1$, generated by a curve
$c(s)=(\rho(s),\lambda(s))$ is given by
\begin{equation}\label{paramHrotaesfgeod}
\left\{\begin{array}{lll}
\!\!\!\sinh^{2}\rho(s) = \displaystyle\frac{A+B \sin(2 \alpha s)}{2 \alpha^{2}},\\
~\\
 \!\!\!\lambda(s)\!=\!\!\displaystyle\int_{0}^{s}\!\!\frac{\sqrt{2} \alpha
(2a+H (1+B \sin(2 \alpha t)))
\sqrt{A+B \sin(2 \alpha t)}}{(2a+H (1+B \sin(2 \alpha t)))^{2}+ \alpha^{2} B^{2} \cos^{2}(2 \alpha t)} ~ dt,\\
~\\
\!\!\!\varphi(t)=t,
\end{array}
\right.
\end{equation}
where $A=1+2aH, ~ B=\sqrt{1+4aH+4 a^{2}}$ and
$\alpha=\sqrt{H^{2}-1}.$

Next we present some notations and definitions we will use
throughout this section. From \eqref{cartescilin} the profile
curve of a rotational $H$-surface in $\mathbb{R}^{3}_{+}$ is given
by
\begin{equation}\label{paramgeratr}
    c_{+}(s) = e^{\lambda(s)}(\tanh \rho (s) ,  \textrm{sech}
   \rho (s)).
\end{equation}

If $H = 1$,~ $\rho(s)$ and $\lambda(s)$ are given by
\eqref{paramH1} and $a>-\displaystyle\frac{1}{2}$. When
$-\displaystyle\frac{1}{2} < a < 0$ we call the family of such of
rotational as  {\it catenoid's  cousin type surfaces} and when
$a=0$ we have the umbilical surfaces with $H = 1$.

If $0 \leq H <1$,~ $\rho(s)$ and $\lambda(s)$ are given by
\eqref{paramHrotaequid} and $a \in \mathbb R$. When $a < 0$  such
rotational surfaces are called {\it equidistant type surfaces} and
when  $a=0$ we get  umbilical surfaces, with $0 \leq H <1$.

If $H > 1$,~ $\rho(s)$ and $\lambda(s)$ are given by
\eqref{paramHrotaesfgeod} and  $a \geq
\displaystyle\frac{-H+\sqrt{H^{2}-1}}{2}$. When
$\displaystyle\frac{-1}{4 H} < a < 0$ we obtain the so called
 {\it onduloid type surfaces} and if $a=0$ we get the umbilical surfaces with $H > 1$.

We observe that the profile curve $c_{+}(s)$ depends on the
parameters $H$, $a$ and $r$, where $r$ is the geodesic radius,
that is,  the radius of a geodesic perpendicular to the
$\overline{z}$-axis.  Then  for $H$ and $r$ fixed we get a family
of rotational $H$-surfaces. From the Euclidean homothety
$\mathcal{H}_{r}$, with factor $r$,
\begin{displaymath}
  \mathcal{H}_{r}(c_{+}(s)) = e^{\lambda(s)}(r \tanh \rho (s) , r \textrm{sech}
   \rho (s)),
\end{displaymath}
we get also other families of profile curves of rotational
$H$-surfaces, so that $r=1$.

According to Theorem 3.1 and Theorem \ref{teo:Alex}  the boundary
of the isoperim\-etric solutions must be rotational $H$-surfaces
that meet the horospheres $\{z=c_{1}\}$ and $\{z=c_{2}\}$
perpendicularly. Our goal is to determine the vertical tangency
points of the profile curves of the rotational surfaces.

\begin{definition}\label{def:tang} Let $c_{+}(s)$ be a curve parametrized by
\eqref{paramgeratr}. We say that a point $c_{+}(s)$  is a vertical
tangency point if the tangent vector in $c_{+}(s)$ satisfies
$\dot{c}_{+}(s)=(0,b)$, where $b \in \mathbb{R}^{\ast}$, that is,

\begin{equation}\label{tgvert}
    e^{\lambda(s)}(\tanh \rho (s) \dot{\lambda}(s)+ \textrm{\textrm{sech}}^{2} \rho(s)
    \dot{\rho}(s))=0,
\end{equation}
\begin{equation}\label{tgvert2}
    e^{\lambda(s)}( \textrm{sech} \rho (s) \dot{\lambda}(s) - \textrm{\textrm{sech}}\rho(s) \tanh \rho (s) \dot{\rho}(s))=b.
\end{equation}
\end{definition}

As  $e^{\lambda(s)} > 0$, \eqref{tgvert} implies that
\begin{equation}\label{tgvert1}
  \tanh \rho (s) \dot{\lambda}(s)+ \textrm{sech}^{2} \rho(s)
  \dot{\rho}(s)=0.
\end{equation}

By \eqref{tgvert1} we obtain the points where the vertical
tangency occurs and by \eqref{tgvert2} we get the direction of the
vertical tangency (upward or downward).

Replacing  \eqref{coordnat} in \eqref{tgvert1} we see that if $p$
is a vertical tangency point for which $U(s) \neq 0$ then
\begin{equation}\label{U}
  U^{2}(s) = \dot{U}^{2}(s),
\end{equation}
and the roots of \eqref{U} give us the vertical tangency points.

Next we study the behaviour of the profile curve of rotational
$H$-surfaces determining the possible vertical tangency points.

\vspace{0.5cm}

{\textbf{First Case: $H=1$}.}

The geodesic radius is given by $\{ \lambda =0\}$ in cylindrical
coordinates. By \eqref{paramgeratr} we get the curve $c_{g}(s) =
(\tanh \rho (s) , \textrm{sech}
   \rho (s))$.

\begin{theorem}\label{teo:simH1}
If $c_{+}(s) = e^{\lambda(s)}(\tanh \rho (s) ,  \textrm{sech}
 \rho (s))$ with $\rho(s), \lambda(s)$ given by \eqref{paramH1} is the parametrization of
 the profile curve of a rotational $H$-surface in $\mathbb{R}^{3}_{+}$,
 with $H=1$, then $c_{+}(s)$ is symmetric to the geodesic radius
 $c_{g}$.
 \end{theorem}
\textbf{ Proof}:~ By \eqref{paramH1} we have that $\lambda(0) =
0$. So $c_{+}(0) \in c_{g}$. If $I$ denotes the Euclidean
inversion through $c_{g}$, we have that $ I(c_{+}(s)) =
c_{+}(-s)$, since $\rho(s)$ is an even function and $\lambda(s)$
is odd.
 \fimprova

 By \eqref{paramH1} we have that
\begin{equation}\label{sinh0}
 \sinh \rho(s) = 0 \Longleftrightarrow a=0 ~ \textrm{e} ~ s=0.
\end{equation}

So $\tanh \rho (s) >0$ if $a$ and $s$ were not both equal to zero.
Furthermore $s=0$ is the unique minimum point of $\rho (s)$.

If $H=1$ we have from \eqref{paramH1} that
\begin{equation}\label{Uhoro}
  U^{2}(s) =
 \displaystyle\frac{a^{2}+(1+2a)^{2}s^{2}}{1+2a}, ~~~ \dot{U}^{2}(s) =
 \displaystyle\frac{(1+2a)^{3}s^{2}}{a^{2}+(1+2a)^{2}s^{2}}.
\end{equation}

Replacing \eqref{Uhoro} in \eqref{U} we get
\begin{equation}\label{eqhoro}
(1+2a)^{4}s^{4} + \Big( 2 a^{2} (1+2a)^{2} - (1+2a)^{4}  \Big)
s^{2} +a^{4} = 0.
\end{equation}

Making $ t = s^{2}$ in \eqref{eqhoro} we get a second degree
equation whose discriminant is
\begin{equation}\label{deltahoro}
\Delta = (1+2a)^{6} (4a +1).
\end{equation}

Since  $1 + 2a > 0$ in case $H=1$ we have that
\begin{itemize}
    \item if $-\displaystyle\frac{1}{2} < a < -\displaystyle\frac{1}{4}$, \eqref{eqhoro} has no real roots. Then
    there are no vertical tangency points in this case;
    \item if  $a = -\displaystyle\frac{1}{4}$, there are at most two vertical tangency points
              \begin{equation}\label{raizeslimhoro}
               s= \pm \displaystyle\frac{1}{2};
              \end{equation}
    \item if $a > -\displaystyle\frac{1}{4}$, there are at most four vertical tangency
    points given by
            \begin{equation}\label{raizeshoro}
               \begin{array}{ll}
               s_{1} =\displaystyle\frac{1+2a+\sqrt{1+4a}}{2 (1+2a)}, ~
               s_{2}=-s_{1},\\
               ~\\
               s_{3} =\displaystyle\frac{1+2a-\sqrt{1+4a}}{2 (1+2a)}, ~
               s_{4}=-s_{3}.
               \end{array}
              \end{equation}

\end{itemize}

Besides these informations we  study the vertical tangencies
according to the variation of the parameter $a$.

\begin{enumerate}
    \item If $-\displaystyle\frac{1}{4} \leq a < 0$ we have
    $\dot{\lambda}(s) >0$. If  $s \geq 0$ then
    \begin{displaymath}
    \tanh \rho (s) \dot{\lambda}(s)+ \textrm{sech}^{2} \rho(s)
  \dot{\rho}(s)>0
  \end{displaymath}
   and  \eqref{tgvert1} is not possible. As $1 +2a
  >0$ the roots $s_{1}$ and $s_{3}$ of \eqref{U} given by  \eqref{raizeshoro} are strictly positives. So they do not give
  vertical tangency points. The other roots $s_{2}, ~ s_{4} < 0$ give us the vertical tangency points that point out
  upward for $b > 0$ in \eqref{tgvert2}. In Fig\-ure \ref{Fig:geraH1aneg} we have the profile curve for $H=1$ and $a =
  -0.2$ and the horocycles that pass through the vertical tangency. In Figure \ref{Fig:rotH1aneg} we see
  two parallel horospheres and the rotational surface, between the horospheres, that meets them perpendicularly.

\begin{figure}[htb] 
       \centering
       \begin{minipage}[b]{0.5\linewidth}
       \centering
          \includegraphics[width=0.8\linewidth]{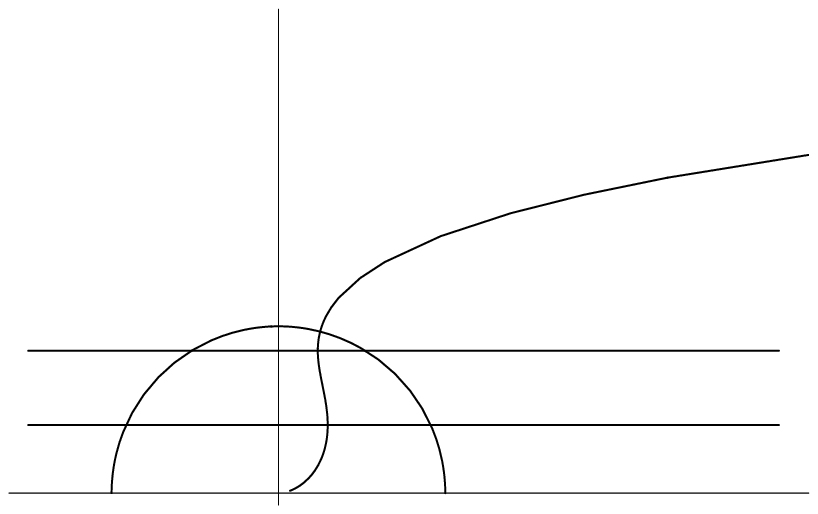}
           \caption{\it Profile curve for $H=1$ and $a = -0.2$.}
           \label{Fig:geraH1aneg}
       \end{minipage}%
              \begin{minipage}[b]{0.5\linewidth}
       \centering
          \includegraphics[width=0.8\linewidth]{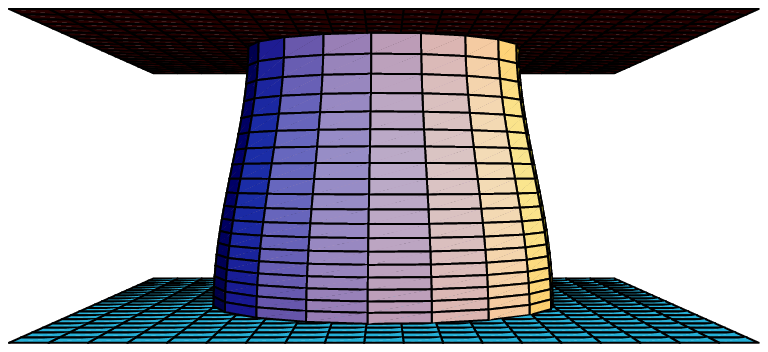}
           \caption{\it Rotational surface with $H=1$ and $a = -0.2$.}
           \label{Fig:rotH1aneg}
       \end{minipage}
   \end{figure}

In particular if  $a = -\displaystyle\frac{1}{4}$, the positive
root $s =
  \displaystyle\frac{1}{2}$ of \eqref{U} given by  \eqref{raizeslimhoro} does not give a vertical tangency point. There exists only
  one vertical tangency point corresponding to $s = -
  \displaystyle\frac{1}{2}$. From the informations about
  $\rho(s)$ and $\lambda(s)$ we get

\begin{equation}\label{comportagerhoro}
 \begin{array}{ll}
\displaystyle\lim_{s \rightarrow -\infty } e^{\lambda(s)} \tanh
\rho (s)  = 0, ~ \displaystyle\lim_{s \rightarrow -\infty }
e^{\lambda(s)} \textrm{sech} \rho (s)  = 0,\\
~\\
\displaystyle\lim_{s \rightarrow \infty } e^{\lambda(s)} \tanh
\rho (s)  = \infty, ~ \displaystyle\lim_{s \rightarrow -\infty }
e^{\lambda(s)} \textrm{sech} \rho (s)  = \infty.
\end{array}
\end{equation}
Therefore this case is not suitable for the isoperimetric problem.

~\\

     \item If $a=0$ we have two pieces of horocycles  tangent at $(0,1)$ that generate the umbilical surfaces
     with $H=1$. They are represented by the Euclidean plane $\{z=1\}$ or the Euclidean sphere
     with radius $\displaystyle\frac{1}{2}$, tangent to $ \partial \mathbb{R}^{3}_{+}$
     at $(0,0,0)$. In the last case it occurs only one vertical tangency point and in this case the surface  meets only one of the horospheres
     (boundary) perpendicularly. In fact, taking the upper
     Euclidean half sphere that represents the horosphere we get
     the possible isoperimetric solution for the umbilical case
     with $H=1$. In Figure  \ref{Fig:geraumbH1} we see the pair
     of profile curves for umbilical surfaces with $H =1$ and the horocycle that pass through the vertical tangency
     point. Figure \ref{Fig:rotumbH1} illustrates the
     possible isoperimetric solution for the umbilical case with
     $H=1$.

\begin{figure}[htb] 
       \centering
       \begin{minipage}[b]{0.5\linewidth}
       \centering
          \includegraphics[width=0.8\linewidth]{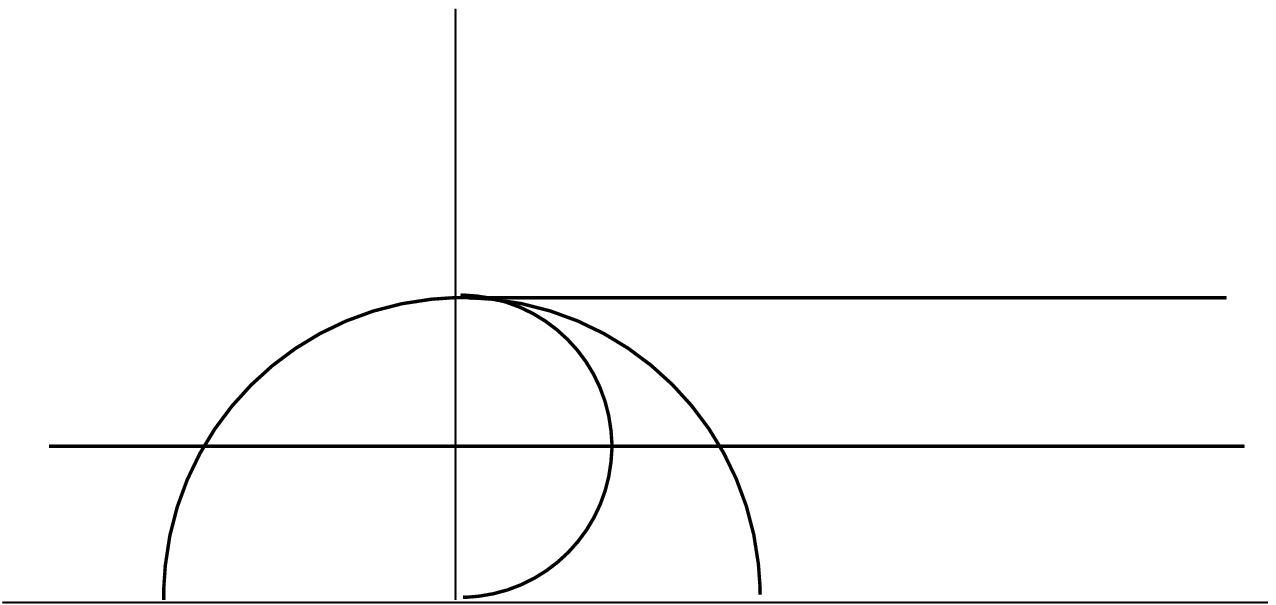}
           \caption{\it Profile
curve for $H=1$ and $a=0$.}
           \label{Fig:geraumbH1}
       \end{minipage}%
              \begin{minipage}[b]{0.5\linewidth}
       \centering
         \includegraphics[width=0.8\linewidth]{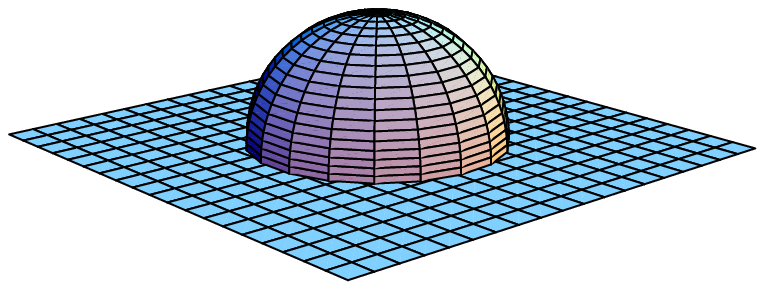}
           \caption{\it Rotational surface with $H=1$ and $a=0$.}
           \label{Fig:rotumbH1}
       \end{minipage}
   \end{figure}

    \item If $a>0$ the profile curves have only one self-intersection. From  \eqref{paramgeratr},  if
    $c_{+}(s_{i}) = c_{+}(s_{j})$ is a self-intersection then
    $s_{i} = \pm s_{j}$ and since the
    curves are symmetric with respect to $c_{g}$, the
    self-intersections must occur on $c_{g}$. So $\lambda (s_{i}) = \lambda (s_{j})
    =0$. From  \eqref{paramH1} we deduce that $\lambda(s)$ has
    a maximum point in $-\displaystyle\frac{\sqrt{a (1+a)}}{1+2a}$
    and a minimum point in $\displaystyle\frac{\sqrt{a
    (1+a)}}{1+2a}$. Furthermore $\displaystyle\lim_{ s \rightarrow \infty} \lambda(s) =
    \infty$ (see \cite {B}), $\lambda (0) =0$ and $\lambda (s)$ is an odd function.

We also have that $\rho(s)$ has only one minimum point in $s=0$.
So if $s > \displaystyle\frac{\sqrt{a (1+a)}}{1+2a}
> 0$ then  $\dot{\rho} (s) , \dot{\lambda }(s) >0$ and
\begin{displaymath}
    \tanh \rho (s) \dot{\lambda}(s)+ \textrm{sech}^{2} \rho(s)
  \dot{\rho}(s)>0;
  \end{displaymath}
if $-\displaystyle\frac{\sqrt{a (1+a)}}{1+2a} < s
 <0$ we have $\dot{\rho} (s) , \dot{\lambda }(s) <0$, which implies that
\begin{displaymath}
    \tanh \rho (s) \dot{\lambda}(s)+ \textrm{sech}^{2} \rho(s)
  \dot{\rho}(s)<0.
  \end{displaymath}
In both cases \eqref{tgvert1} is not verified.  As $a
>0$ the roots $s_{1}, ~ s_{2}, ~ s_{3}, ~ s_{4}$ of
 \eqref{U} given by  \eqref{raizeshoro} satisfy
  \begin{displaymath}
 \begin{array}{ll}
 s_{1} > \displaystyle\frac{\sqrt{a (1+a)}}{1+2a}, ~ s_{2} < - \displaystyle\frac{\sqrt{a
 (1+a)}}{1+2a},\\
 ~\\
 0 < s_{3} < \displaystyle\frac{\sqrt{a (1+a)}}{1+2a}, ~ -\displaystyle\frac{\sqrt{a
 (1+a)}}{1+2a}< s_{4} < 0.
 \end{array}
 \end{displaymath}
 Therefore the vertical tangency is possible only for the positive roots $s_{2}$ and $s_{3}$. As $\dot{\rho} (s_{2}) < 0 $
 and $\dot{\lambda }(s_{2}) > 0$, the vertical tangency in $s_{2}$ points out upward. However $\dot{\rho} (s_{3}) > 0 $ and
  $\dot{\lambda }(s_{3}) < 0$, which implies that the vertical tangency in $s_{3}$ points out downward.
  The isoperimetric solution is not possible in this case because
  if  the vertical tangencies have not occurred in the same height,
  then a piece of the rotational surface would be out of the region
  between the horospheres (see Figure  \ref{Fig:geraH1pos} and \ref{Fig:rotH1pos}).

\begin{figure}[htb] 
       \centering
       \begin{minipage}[b]{0.5\linewidth}
       \centering
           \includegraphics[width=0.8\linewidth]{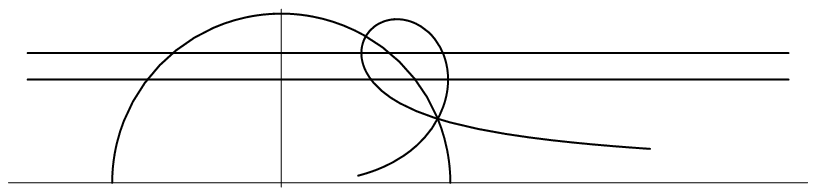}
           \caption{\it Profile
curve for $H=1$ and $a=1$.}
           \label{Fig:geraH1pos}
       \end{minipage}%
              \begin{minipage}[b]{0.5\linewidth}
       \centering
           \includegraphics[width=0.8\linewidth]{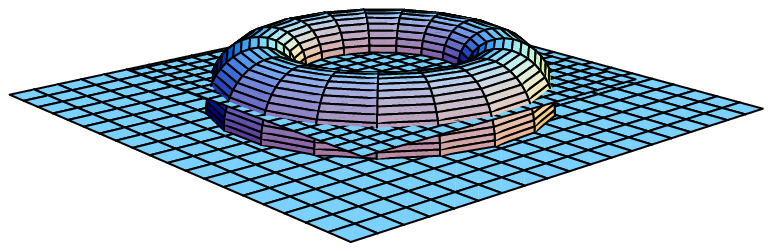}
           \caption{\it Rotational surface with $H=1$ and $a=1$ (excluded).}
           \label{Fig:rotH1pos}
       \end{minipage}
   \end{figure}

Even if the vertical tangency occurred in the same height, the
intersection of the rotational $H$-surface with the parallel
horospheres would be two concentric circles which is not possible
due to Theorem \ref{teo:Alex}.

\end{enumerate}

Then for $H=1$ the boundary of the region $\Omega$ may be a
catenoid's cousin type surface (see Figure  \ref{Fig:rotH1aneg})
or umbilical with $H=1$ (see Figure \ref{Fig:rotumbH1}).\\

 We proceed in the analogous way to study the other cases. We
give here only the main equations and results for them.

\vspace{0.5cm}

{\textbf{Second Case: $0 \leq H <1$}.}

\begin{theorem}\label{teo:simH<1}
If $c_{+}(s) = e^{\lambda(s)}(\tanh \rho (s) ,  \textrm{sech}
 \rho (s))$, with $\rho(s), \lambda(s)$ given by \eqref{paramHrotaequid}, is the parametrization of the profile curve of a rotational $H$-surface in $\mathbb{R}^{3}_{+}$ with $0
\leq H <1$ then
\begin{enumerate}
    \item $c_{+}(s)$ is symmetric with respect to the geodesic radius $c_{g}$;
    \item the assintotic boundary of the profile curves consists of one or two points.
\end{enumerate}
\end{theorem}
\textbf{ Proof}:~ The proof of (a) is similar Theorem
\ref{teo:simH1}. In \cite{B} it is shown that
 $\rho(s)$ is ilimited but $\lambda(s)$ is limited and has finite limit. Then
\begin{displaymath}
\displaystyle\lim_{|s| \rightarrow \infty} e^{\lambda(s)}
\textrm{sech}
 \rho (s) =0,
\end{displaymath}
and the assintotic boundary of the profile curves consists of one
or two points. \fimprova

 As $\sinh \rho(s) \geq 0$, we have from
\eqref{paramHrotaequid} that
\begin{displaymath}
 \sinh \rho(s) = 0 \Longleftrightarrow a=0 ~ \textrm{e} ~ s=0.
\end{displaymath}

If $0 \leq H <1$, from \eqref{paramHrotaequid} it follows that
\begin{equation}\label{Uequid}
  U^{2}(s) =
 \displaystyle\frac{-A+B\cosh(2 \alpha s)}{2 \alpha^{2}}, ~~~ \dot{U}^{2}(s) =
 \displaystyle\frac{ B^{2} \sinh^{2}(2 \alpha s)}{2  (-A+B \cosh(2 \alpha s))}.
\end{equation}

Replacing \eqref{Uequid} in \eqref{U} we get
\begin{equation}\label{eqequid}
B^{2} H^{2} \cosh^{2}(2 \alpha s) - 2 A B \cosh(2 \alpha s) +
A^{2} + \alpha^{2} B^{2} = 0.
\end{equation}

Making $ t = \cosh(2 \alpha s)$ in \eqref{eqequid} we get a second
degree equation whose discriminant is
\begin{displaymath}
\Delta = 4 B^{2} (1-H^{2})^{2} (1+4aH).
\end{displaymath}

As  $B> 0$ in case $0 < H < 1$  and $a$ is definided for any real,
we have that
\begin{itemize}
    \item if $a < - \displaystyle\frac{1}{4 H}$, there are no vertical tangency points;
    \item if $a = -\displaystyle\frac{1}{4 H}$, there are at most two vertical tangency points in
              \begin{equation}\label{raizeslimequid}
               s= \pm \displaystyle\frac{1}{2 \alpha} \textrm{arccosh}\Big(\displaystyle\frac{1}{H}
               \Big);
               \end{equation}
    \item if $a > - \displaystyle\frac{1}{4 H}$, there are at most four vertical tangency points
            \begin{equation}\label{raizesequid}
               \begin{array}{ll}
               s_{1} =\displaystyle\frac{1}{2 \alpha} \textrm{arccosh}\Big(\displaystyle\frac{A+(1-H^{2})\sqrt{1+4 a H}}{B H^{2}}
               \Big), ~
               s_{2}=-s_{1},\\
               ~\\
               s_{3} =\displaystyle\frac{1}{2 \alpha} \textrm{arccosh}\Big(\displaystyle\frac{A-(1-H^{2})\sqrt{1+4 a H}}{B H^{2}}
               \Big), ~
               s_{4}=-s_{3}.
               \end{array}
              \end{equation}
In particular if $H=0$ equation \eqref{eqequid} is written as
\begin{displaymath}
2 B \cosh(2 s) - 1 - B^{2} = 0,
\end{displaymath}
whose solutions are
\begin{displaymath}
               s= \pm \displaystyle\frac{1}{2} \textrm{arccosh}\Big(\displaystyle\frac{B^{2}+1}{2 B}
               \Big).
               \end{displaymath}

\end{itemize}

First, let us suppose $0 < H <1$.

\begin{enumerate}
    \item If $-\displaystyle\frac{1}{4 H} \leq a < 0$, we see that only the roots $s_{2}, ~ s_{4} < 0$ give us the
    vertical tangency points pointing out upward. In Figure \ref{Fig:Hmenoraneg2} we see the profile curve for $H=0.5$ and $a =
  -0.25$ and the horocycles that pass through the vertical tangencies. We observe that the mean curvature vector for the part of the rotational
  surface in the interior of the totally geodesic (symmetry plane of the surface) points out to the rotation axis so determining the isoperimetric
  region illustrated in Figure \ref{Fig:rotHmenoraneg2}.

\begin{figure}[htb] 
       \centering
       \begin{minipage}[b]{0.5\linewidth}
       \centering
          \includegraphics[width=0.8\linewidth]{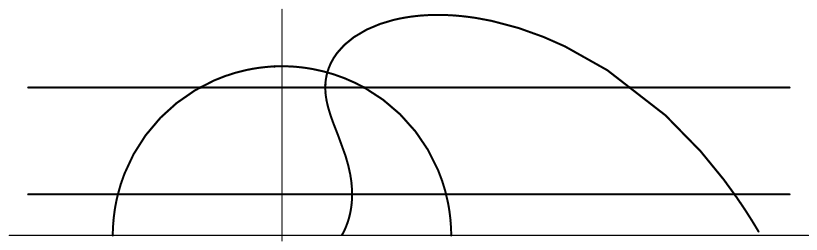}
           \caption{\it Profile curve for $H=0.5$ and $a=-0.25$.}
           \label{Fig:Hmenoraneg2}
       \end{minipage}%
              \begin{minipage}[b]{0.5\linewidth}
       \centering
          \includegraphics[width=0.8\linewidth]{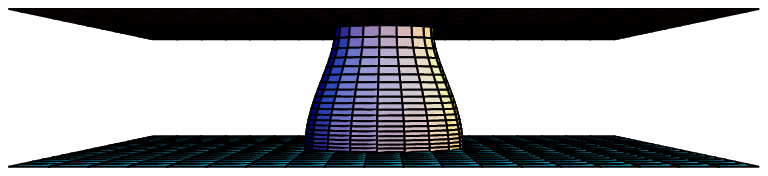}
           \caption{\it Rotational surface with $H=0.5$ and $a=-0.25$.}
           \label{Fig:rotHmenoraneg2}
       \end{minipage}
   \end{figure}

In particular if $a = -\displaystyle\frac{1}{4 H}$ there is only
one vertical tangency point $s = -
  \displaystyle\frac{1}{2 \alpha} \textrm{arccosh}\Big(\displaystyle\frac{1}{H}
  \Big) < 0$. Although the profile curve intersects the horocycle
  in another point, it is not a vertical tangency point.

     \item If $a=0$ we have  two pieces of equidistant curves tangent at $(0,1)$ that generate the umbilical surfaces with
      $0 < H <1$. They are represented by pieces of Euclidean spheres tangent at $(0,0,1)$. We observe that the vertical tangency occurs only for
      the equidistant profile curve nearest the rotation axis. As the mean curvature vector of this umbilical surface points out to the rotation axis
      it determines a isoperimetric region.

    \item If $a>0$ only the roots $s_{2}$ and $s_{3}$ correspond to vertical
tangencies pointing out upward in $s_{2}$ and downward in $s_{3}$,
which implies that this possibility is not suitable for the
isoperimetric problem.
\end{enumerate}

If $H=0$ then for $a<0$ or $a > 0$ we get only one vertical
tangency point which it is not suitable for the isoperimetric
problem (due to the behaviour of the surface). If $a=0$, the
rotational surface is a totally geodesic
    plane which is not suitable for the problem.

Finally we conclude that for $0 \leq H < 1$ the boundary of the
region $\Omega$ may be a equidistant type surface (see Figure
\ref{Fig:rotHmenoraneg2}) or an umbilical surface with $0 < H <
1$.\\

\vspace{0.5cm}

{\textbf{Third Case: $H >1$}.}

\begin{theorem}\label{teo:periH>1}
If $c_{+}(s) = e^{\lambda(s)}(\tanh \rho (s) ,  \textrm{sech}
 \rho (s))$ with $\rho(s), \lambda(s)$ given by \eqref{paramHrotaesfgeod} is the parametrization of the profile
 cuve of a rotational $H$-surface in $\mathbb{R}^{3}_{+}$, with
$H>1$, then $c_{+}(s)$ is a periodic curve with period
$\pi/\alpha$.
\end{theorem}
\textbf{ Proof}:~ We show that the hyperbolic length of the
segment between the points $c_{+}(s)$ and $c_{+}(s+\pi/\alpha)$
keeps constant
 for any $s$. In \cite{B} it was shown that
\begin{equation}\label{rolambda}
\rho(s+\pi/\alpha) = \rho(s),~ \lambda(s+\pi/\alpha)= \lambda(s) +
\lambda(\pi/\alpha),
\end{equation}
which implies from \eqref{paramgeratr} that
\begin{displaymath}
 c_{+}(s+\pi/\alpha)= e^{\lambda(\pi/\alpha} c_{+}(s).
\end{displaymath}

We fix $s_{0}$ and parametrize the segment between the points
$c_{+}(s_{0})$ and $c_{+}(s_{0}+\pi/\alpha)$ by
\begin{displaymath}
\beta(t) = (t, \displaystyle\frac{1}{\sinh \rho (s_{0})} t),~
\textrm{com} ~  e^{\lambda(s_{0})}\tanh \rho (s_{0}) \leq t \leq
e^{\lambda(s_{0})}  e^{\lambda(\frac{\pi}{\alpha})} \tanh \rho
(s_{0}).
\end{displaymath}
So its hyperbolic length is
\begin{displaymath}
    L(\beta(t)) = \lambda(\displaystyle\frac{\pi}{\alpha}) \cosh
    \rho (s_{0}).
\end{displaymath}

The length of the segment depends only on the function $\rho(s)$
with period $\pi/\alpha$, given in \eqref{rolambda}. So
$L(\beta(t))$ is constant for any $s_{0}$. \fimprova

If $H >1$ we have from  \eqref{paramHrotaesfgeod} that
\begin{equation}\label{Uesfgeod}
  U^{2}(s) =
 \displaystyle\frac{A+B\sin(2 \alpha s)}{2 \alpha^{2}}, ~~ \dot{U}^{2}(s) =
 \displaystyle\frac{B^{2} \cos^{2}(2 \alpha s)}{2 (A+B \sin(2 \alpha s))}.
\end{equation}

Replacing \eqref{Uesfgeod} in  \eqref{U} we get
\begin{equation}\label{eqesfgeod}
B^{2} H^{2} \sin^{2}(2 \alpha s) + 2 A B \sin(2 \alpha s) + A^{2}
- \alpha^{2} B^{2} = 0.
\end{equation}

Making $ t = \sin(2 \alpha s)$ in \eqref{eqesfgeod}, the
discriminant of equation \eqref{eqesfgeod} is
\begin{equation}\label{deltaesfgeod} \Delta = 4 B^{2}
(H^{2}-1)^{2} (1+4aH).
\end{equation}

From \eqref{deltaesfgeod} we have the following analysis.

As $H > 1$ and $B> 0$,
\begin{itemize}
    \item if $a < - \displaystyle\frac{1}{4 H}$ there are no vertical tangency points;
    \item if $a = - \displaystyle\frac{1}{4 H}$ the possible vertical
    tangency points are
              \begin{equation}\label{raizeslimesfgeod}
               \begin{array}{ll}
               s_{k}= s_{0} + \displaystyle\frac{k \pi}{\alpha}, ~ k \in
               \mathbb{Z},\\
               ~\\
               \tilde{s}_{k}=\tilde{s}_{0} + \displaystyle\frac{k \pi}{\alpha} , ~
               k \in \mathbb{Z},
               \end{array}
               \end{equation}
where \begin{displaymath} \begin{array}{ll}
s_{0}=\displaystyle\frac{1}{2 \alpha}
\arcsin\Big(\displaystyle\frac{1}{ H} \Big), ~
\displaystyle\frac{- \pi}{2} < 2 \alpha s_{0} <
\displaystyle\frac{ \pi}{2},\\
~\\
\tilde{s}_{0} = \displaystyle\frac{1}{2 \alpha}
\arcsin\Big(\displaystyle\frac{1}{ H} \Big), ~
\displaystyle\frac{\pi}{2} < 2 \alpha \tilde{s}_{0} <
\displaystyle\frac{3 \pi}{2}.
\end{array}
\end{displaymath}
    \item if $a > -\displaystyle\frac{1}{4 H}$ the possibilities are
            \begin{displaymath}
               \begin{array}{llll}
               S_{k} = S_{0} + \displaystyle\frac{k \pi}{\alpha}, ~ k \in
               \mathbb{Z},\\
               ~\\
               \tilde{S}_{k}=\tilde{S}_{0} + \displaystyle\frac{k \pi}{\alpha} , ~
               k \in \mathbb{Z},\\
               ~\\
                s_{k} = s_{0} + \displaystyle\frac{k \pi}{\alpha}, ~ k \in
               \mathbb{Z},\\
               ~\\
               \tilde{s}_{k}=\tilde{s}_{0} + \displaystyle\frac{k \pi}{\alpha} , ~
               k \in \mathbb{Z},
               \end{array}
              \end{displaymath}
where
\begin{displaymath}
\begin{array}{llll}
S_{0} = \displaystyle\frac{1}{2 \alpha}
\arcsin\Big(\displaystyle\frac{-A+(H^{2}-1) \sqrt{1+4aH}}{B H^{2}}
\Big), ~ \displaystyle\frac{- \pi}{2} < 2 \alpha S_{0} <
\displaystyle\frac{\pi}{2},\\
~\\
\tilde{S}_{0} = \displaystyle\frac{1}{2 \alpha}
\arcsin\Big(\displaystyle\frac{-A+(H^{2}-1) \sqrt{1+4aH}}{B H^{2}}
\Big), ~ \displaystyle\frac{\pi}{2} < 2 \alpha \tilde{S}_{0}
< \displaystyle\frac{3 \pi}{2},\\
~\\
s_{0} = \displaystyle\frac{1}{2 \alpha}
\arcsin\Big(\displaystyle\frac{-A-(H^{2}-1) \sqrt{1+4aH}}{B H^{2}}
\Big), ~ \displaystyle\frac{- \pi}{2} < 2 \alpha s_{0} <
\displaystyle\frac{\pi}{2},\\
~\\
\tilde{s}_{0} = \displaystyle\frac{1}{2 \alpha}
\arcsin\Big(\displaystyle\frac{-A-(H^{2}-1) \sqrt{1+4aH}}{B H^{2}}
\Big)$, ~ $\displaystyle\frac{\pi}{2} < 2 \alpha \tilde{s}_{0} <
\displaystyle\frac{3 \pi}{2}.
\end{array}
\end{displaymath}
\end{itemize}

Now we determine when the vertical tangency really occurs,
depending on the geometry of the profile curve.

\begin{enumerate}
    \item If $-\displaystyle\frac{1}{4 H} \leq a < 0$  only the roots $\tilde{S}_{k}, ~ \tilde{s}_{k}$
furnish the vertical tangency points pointing out upward. In
Figure \ref{Fig:onduloide} we see a hyperbolic onduloid.

\begin{figure}[ht]
\begin{center}
\includegraphics[width=2cm]{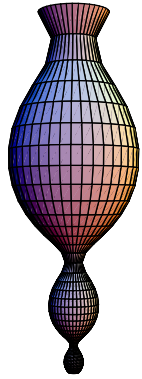}
\caption{ Hyperbolic  onduloid with $H =3$ and $a = -0.05$.}
\label{Fig:onduloide}
\end{center}
\end{figure}

     \item If $a=0$ we have tangent geodesic half circles along de rotation
     axis. Each geodesic half circle generates a geodesic sphere in
     $\mathbb{R}^{3}_{+}$ which are  umbilical surfaces with $H > 1$ and isoperimetric regions.

    \item If $a>0$ we analyse the behaviour of
the profile curve in the interval $\Big] - \pi/4 \alpha, 3 \pi/4
 \alpha\Big[$ since according to Theorem \ref{teo:periH>1} the profile curves are periodic with
period $\pi/\alpha$. It is easy to see that only the roots  $
s_{0}$ and $\tilde{S}_{0}$ furnish the vertical tangency pointing
out downward in $ s_{0}$ and upward in $\tilde{S}_{0}$.
  This case corresponds to the so called hyperbolic nodoids, which
 are not embedded surfaces.
\end{enumerate}

We conclude that for $H > 1$ the boundary of the region $\Omega$
can be an onduloid's type surface (see Figure
 \ref{Fig:onduloide}) or an umbilical surface with $H > 1$.

Next we prove  Theorem 1.1. We start with the existence an then we
obtain the possible minimizing regions.

\textbf{Proof of Theorem 1.1}:~ By Theorem \ref{teo:Alex} , the
solutions to the isoperimetric problem have as boundaries
rotationally invariant surfaces which have constant mean curvature
where they are regular. But they must be regular (actually
analytic), since the singularities along such boundaries must
have, by well-known results, (Hausdorff) codimension at least 7,
which is not possible for (2-dimensional) surfaces. Now, from
\cite{Mo1} the existence of the isoperimetric solutions follows
from the fact that $\mathcal{F}_{c_{1}, c_{2}}/G$ is compact,
where $G$ is the group of isometries of $\mathbb{R}^{3}_{+}$ whose
elements let invariant the region $\mathcal{F}_{c_{1}, c_{2}}$
between the horospheres, that is, the rotations around a vertical
geodesic and the horizontal translations. The second part of
Theorem 1.1 follows from the analysis of vertical tangencies done
in the First, Second and Third Cases above.


\end{document}